\documentclass[11pt,a4paper]{article}

\title{Generalized stochastic Lagrangian paths for the Navier-Stokes equation}
\author{Marc Arnaudon$^a$\quad Ana Bela Cruzeiro$^b$\quad and Shizan Fang$^c$
\vspace{3mm}\\
{\footnotesize $^a$Institut de Maths. de Bordeaux,
Universit\'e de Bordeaux I,
33405 Talence Cedex, France}\\
{\footnotesize $^b$GFMUL and Dep. de Mat.  Instituto Superior T\'ecnico, Av. Rovisco Pais,
1049-001 Lisboa, Portugal }\\
{\footnotesize $^c$Institut de Maths. de Bourgogne,
Universit\'e de Bourgogne,
21078 Dijon Cedex, France}}
\date{February 18th,  2016}

\usepackage{amssymb,amsmath,amsfonts,amsthm,color}

\usepackage{framed}

\theoremstyle{plain}
\newtheorem{thm}{Theorem}[section]

\newtheorem{prop}[thm]{Proposition}

\theoremstyle{definition}
\newtheorem{definition}[thm]{Definition}
\newtheorem{remark}[thm]{Remark}

\newtheorem{example}[thm]{Example}

\catcode`\@=11
\def\mequal{\mathrel{\mathpalette\@mvereq{\hbox{\sevenrm m}}}} 
\def\@mvereq#1#2{\lower.5\p@\vbox{\baselineskip\z@skip\lineskip1.5\p@
    \ialign{$\m@th#1\hfil##\hfil$\crcr#2\crcr=\crcr}}}
\def\partr#1#2{/\kern-.08333em/_{#1,#2}^{\phantom{.}}}
\def\invpartr#1#2{/\kern-.08333em/_{#1,#2}^{-1}} 
\def\hpartr#1#2{/\kern-.08333em/_{#1,#2}^{h}}
\def\Epartr#1#2{/\kern-.08333em/_{#1,#2}^{E}}

\def\newdot{{\kern.8pt\cdot\kern.8pt}}
\def\,{\relax\ifmmode\mskip\thinmuskip\else\thinspace\fi}
\def\{{\relax\ifmmode\lbrace\else $\lbrace$\fi}
\def\}{\relax\ifmmode\rbrace\else $\rbrace$\fi}
\def\R{\mathbb{R}}
\def\L{\mathcal L}
\def\E{\mathbb{E}}
\def\Z{\mathbb{Z}}
\def\T{\mathbb{T}}
\def\F{\mathcal{F}}
\def\Ric{\textup{Ric}}
\def\div{\textup{div}}
\def\Diff{{\rm Diff}}
\def\dis{\displaystyle}

\newcommand{\ra}{\rightarrow}
\newcommand{\ee}{\varepsilon}

\font\sevenrm=cmr7

\renewcommand{\theequation}{\arabic{section}.\arabic{equation}}

\begin{document}
\maketitle
\makeatletter 
\renewcommand\theequation{\thesection.\arabic{equation}}
\@addtoreset{equation}{section}
\makeatother 



%
%

\begin{abstract}

\it  In the note added in proof of the seminal paper [Groups of diffeomorphisms and the motion of an incompressible fluid,
{\it Ann. of Math.} 92 (1970), 102-163], Ebin and Marsden introduced the so-called {\it correct Laplacian} for the Navier-Stokes equation on 
a compact Riemannian manifold. In spirit of Brenier's generalized flows for the Euler equation, 
 we  introduce a class of semimartingales on a compact Riemannian manifold. We prove that these semimartingales are critical points to  the corresponding kinetic energy
 if and only if its drift term solves weakly the Navier-Stokes equation defined with Ebin-Marsden's Laplacian. 
We also show that for the case of torus,  classical solutions of the Navier-Stokes equation realize the minimum of the kinetic energy in a suitable class.\rm

\end{abstract}

\section{ \bf Introduction}\label{section1}

\quad Euler equations describe the velocity of incompressible non-viscous fluids. Considering these equations  on a bounded domain $U$ of $\R^d$, or 
on a compact Riemannian manifold $M$ without boundary, they read

\begin{equation}\label{eq1.1}
{d\over dt}u_t+ (u_t\cdot\nabla)u_t=-\nabla p,\quad \div(u_t)=0.
\end{equation}

Lagrange's point of view consists in  describing the position of the particles: for a solution  $u$ to \eqref{eq1.1},
 it concerns solutions of the ordinary differential equation (ODE)
\begin{equation}\label{eq1.2}
{d\over dt}g_t(x)=u_t(g_t(x)),\quad g_0(x)=x.
\end{equation}

When $(t,x) \ra u_t(x)$ is smooth, the ODE \eqref{eq1.2} defines a flow of $C^\infty$-diffeomorphisms $g_t$. From the position values, we get the velocity by
$$\dis u_t(x)=\big( {d\over dt}g_t \big) (g_t^{-1}(x)).$$

In this case, the two points of view are equivalent. 
Throughout the whole paper  we shall consider the interval of time $[0,T]$.  Equation \eqref{eq1.2} defines a continuous map
$$ g_\cdot:  [0,T]\ra \Diff(M)$$
from $[0,T]$ to the group of diffeomorphisms of $M$.
\vskip 2mm
In a famous work \cite{Arnold}, V.I.  Arnold gave a geometric interpretation to the incompressible Euler equation, saying that  $u$ is a solution to \eqref{eq1.1}
if and only if $t\ra g_t$ is a {\it  geodesic} on the submanifold of $\Diff(M)$ keeping the volume measure invariant, 
equipped with the $L^2$ metric.  Equivalently, $g_\cdot$ minimizes the action

\begin{equation}\label{eq1.3}
S[\varphi]={1\over2}\int_0^T\!\!\int_M \Bigl|{d\over dt}\varphi_t(x)\Bigr|_{T_xM}^2\,dxdt
\end{equation}
on $C([0,T], \Diff(M))$, where $dx$ denotes the normalized Lebesgue measure on $U$ or the normalized Riemannian volume on $M$ (see also \cite{EM}).

\vskip 2mm
In \cite{Brenier}, Y. Brenier gave a probabilistic interpretation to \eqref{eq1.1}, by  looking for probability
measures $\eta$ on the path space $C([0,T],M)$, which minimize the kinetic energy

\begin{equation}\label{eq1.4}
S[\eta]={1\over 2}\int_{C([0,T],M)}\Bigl[\int_0^T |\dot \gamma(t)|_{T_{\gamma(t)}M}^2\,dt\Bigr]d\eta(\gamma),
\end{equation}
with constraints $(e_t)_*\eta=dx$, where  $e_t: \gamma\ra \gamma(t)$ denotes the evaluation map. Let 
$$X(\gamma, t)=\gamma(t).$$
Then under $\eta$, $\{X(\cdot, t); t\geq 0\}$ is a $M$-valued stochastic process. Moreover, in \cite{Brenier} as well as in \cite{Brenier2}, 
Brenier proved that such a probability measure $\eta$ gives rise 
to a weak solution of the Euler equation in the sense of Di Perna and Majda  \cite{DiPerna}.  More precisely, define
a probability measure $\mu$ on $[0,T]\times TM$ by
\begin{equation*}
\begin{split}
&\int_{[0,T]\times TM} f(t,x,v)\,\mu(dt,dx,dv)\\
&={1\over T}\int_0^T\!\!\int_{C([0,T],M)} f(t, \gamma(t), \gamma'(t))d\eta(\gamma)\,dt.
\end{split}
\end{equation*}
Then $\mu$ solves the Euler equation in generalized sense:
\begin{equation*}
\int \Bigl[ v\cdot w(x)\alpha'(t) + v\cdot (\nabla w(x)\cdot v)\alpha(t)\Bigr]\, \mu(dt,dx,dv)=0
\end{equation*}
for any $\alpha\in C^\infty(]0,T[)$ and any smooth vector field $w$ such that $\div(w)=0$.  We also refer to \cite{AmbrosioF} in which the authors used 
the theory of  mass transportation.

\vskip 2mm

In this work, we will  deal with Navier-Stokes equations on a compact Riemannian manifold $M$. There are two natural ways to define the ``Laplace'' operator 
on vector fields. The first way is to use the de Rham-Hodge Laplace operator $\square$ on differential $1$-forms, that is 
$\square=dd^*+d^*d$. As usual, for a vector field $A$, we denote by
$A^\sharp$ the associated differential $1$-form; for a differential $1$-form $\omega$, we denote by $\omega^\flat$ the corresponding vector field. Then we define
\begin{equation*}
\square A=(\square A^\sharp)^\flat.
\end{equation*}
The Weitzenb\"ock formula states that 
\begin{equation}\label{eq1.4-1}
-\square A=\Delta A- \Ric (A)
\end{equation}
where $\Delta A= \hbox{\rm Trace}(\nabla^2 A)$ and $\Ric$ is the Ricci tensor. Another natural way, following \cite{EM},
 is to use the deformation tensor. More precisely, let $A$ be a vector field on $M$,
the deformation tensor $\hbox{\rm Def}\,A$ is a symmetric tensor of type $(0,2)$ such that
\begin{equation}\label{eq1.4-2}
(\hbox{\rm Def}\,A)(X,Y)={1\over 2}\Bigl( \langle\nabla_X A, Y\rangle + \langle\nabla_Y A,X\rangle\Bigr).
\end{equation}

Then $\hbox{\rm Def}: TM\ra S^2T^*M$ which sends a vector field to a symmetric tensor of type $(0,2)$. Let $\hbox{\rm Def}^*: S^2T^*M\ra TM$ be the adjoint operator.
According to \cite{MT}, as well as to \cite{Taylor}, we define
\begin{equation}\label{eq1.4-3}
\hat\square = 2\hbox{\rm Def}^*\hbox{\rm Def}.
\end{equation}
Then on vector fields of divergence free $A$, it holds true (see \cite{MT, Pier})
\begin{equation}\label{eq1.4-4}
-\hat\square A= \Delta A + \Ric (A).
\end{equation}

Comparing \eqref{eq1.4-4} to \eqref{eq1.4-1}, the sign of $\Ric$ is opposite. 
\vskip 2mm

In this work, we will consider the following Navier-Stokes equation on  $M$
\begin{equation}\label{eq1.5}
{d\over dt}u_t + \nabla_{u_t}u_t +\nu\, \hat\square u_t=-\nabla p,\ \div(u_t)=0,
\end{equation}

where $\nu>0$ is the viscosity coefficient. Since $\div(u_t)=0$, we have $\int_M \langle \nabla_{u_t}u_t, u_t\rangle\, dx=0$.
Using the relation $\hat\square=\square -2\Ric$ and equation \eqref{eq1.5}, we get

\begin{equation}\label{eq1.5-1}
{1\over 2}\, {d\over dt}\int_M |u_t|^2\, dt +\nu\int_M \Bigl(|du_t^\sharp|^2 +|d^*u_t^{\sharp}|^2\Bigr)\, dx 
-2\nu \int_M \langle \Ric\, u_t, u_t\rangle\, dx=0.
\end{equation}

When $\Ric$ is negative, the above relation yields the existence of Leray's weak solution (see for instance \cite{Pier}).  
For the general case, the existence of Leray's weak solution to \eqref{eq1.5} was proved in \cite{Taylor} (Theorem 4.6, p.498 and p.504).
\vskip 2mm

In contrast to Euler equations, there is no geometrical interpretation for Navier-Stokes equations.
The purpose of this work is to develop a probabilistic interpretation to equation \eqref{eq1.5}. Note that in this context, 
it is suitable to consider that the underlying Lagrangian trajectories are  semimartingales $\xi_t$ on the manifold $M$. 
Comparing to Brenier's generalized flows for Euler equations,  the paths $t\ra \xi_t$ 
are never of finite energy in the sense of \eqref{eq1.3}. Instead, we shall consider the mean kinetic energy (see definition \eqref{eq2.7} below). This functional 
first appeared in stochastic optimal control \cite{Fleming} as well as in connection with quantum mechanics \cite{Zambrini}; we mention also  \cite{Follmer} for the relation of (stochastic) kinetic energy and entropy as well as
 \cite{Yasue}, for its appereance in the study of the Navier-Stokes equation. 
 
 \vskip 2mm
 
 Roughly speaking, the main result of this paper (see Theorem \ref{th1} below) says that the semi-martingale $\xi_t$ in a suitable class is a critical  point to the stochastic kinetic energy \eqref{eq2.7}
 if and only if its drift term $u_t$ solves Navier-Stokes equation \eqref{eq1.5} in the sense of Di-Perna and Majda.
 \vskip 2mm
 
  In the recent years  the functional \eqref{eq2.7} has been used with success in various contexts  (see for example \cite{AC1, ACC, AC2, AC3, CCR, Constantin, Cruzeiro, Eyink, Iyer, Leonard}).
  In comparison with \cite{AC1, ACC, AC2, AC3},  we do not require, in the present work,  that martingales have the flow property.

\vskip 2mm
The organisation of the paper is as follows. In section \ref{section2}, we shall introduce and study the class of $\nu$-Brownian incompressible semimartingales. 
We prove that such a semimartingale is a critical point of the corresponding kenetic energy \cite{Cruzeiro} if and only if it solves the Navier-Stokes equation in the 
sense of DiPerna-Majda \cite{Brenier, Brenier2}. We also prove the existence of a minimum under certain conditions. In section \ref{section3}, we shall show, in the case of  a torus $\T^d$,
 that a classical solution to Navier-Stokes equation gives rise to a $\nu$-Brownian incompressible martingale which realizes the minimum of the 
 kinetic energy in a convenient class.
 
 \vskip 2mm

\section{\bf Generalized stochastic paths for the Navier-Stokes equation }\label{section2}

In this section,  $M$ will denote a  connected compact Riemannian manifold without boundary. 
Let $(\Omega, \F, P)$ be a probability space equipped with a filtration $\{\F_t;\, t\geq 0\}$ satisfying the usual conditions.

\vskip 2mm
 A $M$-valued stochastic process $\xi_t$ defined on $(\Omega, \F, P)$ is
said to be a semimartingale on $M$ if for any $f\in C^2(M)$, $f(\xi_t)$ is a real valued semimartingale. This notion is independent of the chosen connection on $M$; however, 
the corresponding  local characteristics are dependent of the choice of connection. For a semimartingale $(\xi_t)$ starting from a point $x\in M$ and given a connection $\nabla$, 
 the stochastic parallel translation $//_t$ along $\xi_\cdot$ can be defined and

$$\dis \zeta_t=\int_0^t //_s^{-1}\circ d\xi_s$$
is a $T_xM$-valued semimartingale. Then  there exist processes $(\xi^0(s), H_1(s), \ldots, H_m(s))$ which are adapted to $\F_t$ such that 
\begin{equation*}
\xi^0(s),H_1(s), \ldots, H_m(s)\in T_{\xi_s}M
\end{equation*}
 and  $\zeta_t$ admits It\^o form
 
\begin{equation}\label{eq2.0}
\zeta_t=\int_0^t  //_s^{-1}\xi^0(s)\, ds + \sum_{i=1}^m \int_0^t  //_s^{-1}H_i(s)\, dw_s^i
\end{equation}
where $w_t=(w_t^1, \cdots, w_t^m)$ is a standard Brownian motion on $\R^m$ (see for example \cite{Bismut}). 
For example, if the semimartingale $\xi_t$ comes from a stochastic differential equation (SDE) on $M$:

\begin{equation*}
d\xi_t=X_0(t,\xi_t) dt + \sum_{i=1}^m X_i(t,\xi_t)\circ dw_t^i,\quad \xi_0=x,
\end{equation*}

then 
\begin{equation*}
\xi^0(t)=X_0(t,\xi_t)+ {1\over 2} \sum_{i=1}^m(\nabla_{X_i}X_i)(t, \xi_t).
\end{equation*}

For simplicity, in what follows, we only consider the Levi-Civita connection $\nabla$ on $M$.
As in  \cite{ACC}, \cite{Cruzeiro}, we consider the operator
\begin{equation}\label{eq2.1}
D_t\xi=//_t\, \lim_{\ee\ra 0}\E\Bigl( {\zeta_{t+\ee}-\zeta_t\over \ee}\Big|\F_t\Bigr),
\end{equation}
\vskip 1mm
which is well-defined and equals  $\xi^0(t)$. For a semimartingale $\xi_t$ given by \eqref{eq2.0}, the It\^o formula has the following form (see \cite{Bismut}, p. 409)

\begin{equation}\label{Ito}
\begin{split}
f(\xi_t)=f(\xi_0)&+\int_0^t \Bigl(\langle \nabla f(\xi_s), \xi^0(s)\rangle +{1\over 2}\sum_{i=1}^m\langle \nabla_{H_i(s)}(\nabla f)(\xi_s), H_i(s)\rangle\Bigr)\,ds\\
&+\sum_{i=1}^m\int_0^t \langle \nabla f(\xi_s), H_i(s)\rangle\, dw_s^i.
\end{split}
\end{equation}
\vskip 2mm

Let $\{g_t (x,\omega );\  t\geq 0, x\in M, \omega \in \Omega \}$ be a family of continuous semimartingales with values in $M$. 
Let $\mathbb P^g$ denote the law of $g$ in 
the continuous path space $C([0,T],M)$, that is, for every cylindrical functional $F$,

\begin{equation*}
\int_{C([0,T], M)} F(\gamma (t_1 ) ,...,\gamma (t_n))d\mathbb P^g (\gamma )
=\int_M \Bigl[ \int _{C([0,T], M)}  F(g_{t_1}(x),..., g_{t_n}(x))d\mathbb P^g_x \Bigr]dx
\end{equation*}

where $\mathbb P^g =\mathbb P^g_x \otimes dx$ and 
under $\mathbb P^g_x$, the semimartingale $g_t$ starts from $x$. 

\vskip 2mm
We shall say that the semimartingale $g_t$ is \it incompressible \rm if, for each $t>0$,

\begin{equation}\label{eq2.2}
\E_{\mathbb P^g}[ f(g_t)]=\int_M f(x)dx, \quad \hbox{for all } f\in C(M)
\end{equation}
the expectation being taken with respect to the law $\mathbb P^g$ of $g$.

\vskip 2mm
Let $\nu>0$; we shall say that $g_t$ is a $\nu$-{\it Brownian semimartingale} if,  under $\mathbb P^{g}$, there
 exists a time-dependent adapted random vector field $u_t$ over $g_t$ such that

\begin{equation}\label{eq2.3}
M_t^f =f(g_t )-f(g_0 )-\int_0^t\Bigl(\nu \Delta f(g_s) + \langle u_s, \nabla f(g_s )\rangle \Bigr)ds,
\end{equation}

is a local continuous martingale with the quadratic variation given by
$$\langle M_t^{f_1} ,M_t^{f_2} \rangle =2\nu \int_0^t \langle \nabla f_1, \nabla f_2\rangle (g_s )ds.$$

For a semimartingale $\xi_t$ given by \eqref{eq2.0},  if $\{H_1(s), \ldots, H_m(s)\}$ is an orthogonal system such that
for any vector $v\in T_{\xi_s}M$, $\dis \sum_{i=1}^m\langle v, H_i(s)\rangle^2 = 2\nu |v|^2$, then it is a $\nu$-{Brownian semimartingale}.

\begin{example}\label{ex1}
In the flat case $\R^d$, such a semimartingale admits the following form
\begin{equation}\label{eq2.4}
dg_t(w)=\sqrt{2\nu}\, dw_t + u_t(w)\, dt,
\end{equation}
where $(w_t)$ is a Brownian motion on $\R^d$ and $\{u_t; t\geq 0\}$ is an adapted $\R^d$-valued process
such that $ \int_0^T |u_t(w)|^2\,dt<+\infty$ almost surely. $\square$
\end{example}

\begin{example}\label{ex2}
For the general case of a compact Riemannian manifold $M$, we consider the bundle of orthonormal frames $O(M)$.
 Let $(V_t)_{t\in [0,T]}$ be a family of $C^1$ vector fields such that  the dependence $t\ra V_t$ is $C^1$. Denote by $\tilde V_t$ the 
 horizontal lift of $V_t$ to $O(M)$. Let $\div (V_t)$ and $\div(\tilde V_t)$ be respectively the divergence operators on $M$ and on $O(M)$;  
 they are linked by (see \cite{FangLuo}, p. 595)
 
 $$\dis \div(\tilde V_t)=\div(V_t)\circ \pi,$$
 where $\pi: O(M)\ra M$ is the canonical projection. It follows that if $\div(V_t)=0$, then $\div(\tilde V_t)=0$. 
 Consider the horizontal diffusion $r_t$ on $O(M)$ defined by the SDE
 
 \begin{equation}\label{diffusion1}
 dr_t=\sqrt{2\nu}\, \sum_{i=1}^d H_i(r_t)\circ dw_t^i + \tilde V_t(r_t)\, dt,\quad r_0\in O(M)
 \end{equation}
 where $\{H_1, \cdots, H_d\}$ are the canonical horizontal vector fields on $O(M)$. 
 Let $dr$ be the Liouville measure on $O(M)$, then the stochastic flow $r_0\ra r_t(r_0)$ leaves $dr$ invariant. 
 Set
 
 \begin{equation}\label{diffusion2}
  \xi(t,x)=\pi(r_t(r_0)),\quad r_0\in \pi^{-1}(x).
 \end{equation}
 
For any continuous function $f$ on $M$,
 $$\dis \int_M\E(f(\xi(t,x))\, dx=\int_M f(x)\, dx.$$
 Then $\xi$ is an incompressible $\nu$-Brownian diffusion, with $D_t\xi(x) =V_t(\xi(t,x)). \ \square$
 \end{example}
 
 \begin{remark} Let $P_t$ be the semigroup associated to $\dis {1\over 2}\Delta_M + V_t$ with $\div(V_t)=0$; then for any $f\in C^2(M)$,
 \begin{equation*}
 {d\over dt}\int_M P_t f(x)\, dx=\int_M \bigl( {1\over 2}\Delta_M P_t f + V_t P_tf\bigr)\, dx=0.
 \end{equation*}
 It follows that for any continuous function $f: M\ra \R$,
 $$\dis \int_M P_tf (x)\, dx=\int_M f(x)\, dx.$$
 Therefore any SDE on $M$ defining a Brownian motion with drift $V$ gives rise to an incompressible $\nu$-Brownian diffusion $\xi$ with $D_t\xi(x) =V_t(\xi(t,x))$.
 
 \end{remark}
  \begin{example}\label{ex3}
 Let $\Z^2$ be the set of two dimensional lattice points and define $\Z^2_0=\Z^2\setminus \{(0,0)^\ast \}$. For $k\in\Z^2_0$, we consider the vector $k^\perp=(k_2,-k_1)^\ast$ and the vector fields
  $$A_k(\theta)=\sqrt{\nu\over\nu_0}\frac{\cos(k\cdot \theta)}{|k|^\beta}k^\perp, \quad B_k(\theta)=\sqrt{\nu\over\nu_0}\frac{\sin(k\cdot \theta)}{|k|^\beta}k^\perp,\quad \theta\in\T^2,$$
where $\beta>1$ is some constant.

Let $\tilde \Z^2_0$ the subset of $\Z^2_0$ where we identify vectors $k,k^\prime $such that $k+k^\prime =0$ and let
$$\nu_0= \sum_{k\in \tilde \Z^2_0}\frac{1}{2|k|^{2\beta}}.$$
The family $\{A_k,B_k:k\in \Z^2_0\}$ constitutes an orthogonal basis of the space of divergence free vector fields on $\T^2$
and satisfies 

\vskip 2mm
$$\dis \sum_{k\in \tilde \Z^2_0}  \Bigl(\langle A_k, v\rangle^2+\langle B_k, v\rangle^2\Bigr)=\nu\, |v|_{T_{\theta}\T^2}^2,\quad v\in T_{\theta}\T^2,$$
and 
$$\dis \sum_{k\in \tilde \Z^2_0}   \nabla_{A_k}A_k  =0,~ \sum_{k\in \tilde \Z^2_0}   \nabla_{B_k}B_k =0.$$
Consider the SDE on $\T^2$,

\begin{equation}\label{A1}
d\xi_t=\sum_{k\in \tilde \Z^2_0}  \Bigl( {A_k}(\xi_t)\circ dw_t^{k} + {B_k}(\xi_t)\circ d\tilde w_t^k\Bigr) + u(t, \xi_t)\,dt,\quad \theta_0=\theta\in\T^2
\end{equation}
where $\{w_t^k, \tilde w_t^k;\ k\in \Z^2_0\}$ are independent standard Brownian motions on $\R$, and $u(t, \cdot)$ is a family of divergence free
vector fields in $H^1(\T^2)$, such that,
$$\dis \int_0^T\!\!\int_{\T^2} ( |u|^2+|\nabla u|^2)\,dxdt<+\infty.$$
Then by  \cite{Cruzeiro, FangLuo}, for $\beta\geq 3$, 
the SDE \eqref{A1} defines a stochastic flow of measurable maps which preserves the Haar measure $dx$ on $\T^2$. More precisely, for almost surely $w$, the map
$$\dis x\ra \xi_t(x,w)\ \hbox{\rm solution to }\eqref{A1} \hbox{ with initial condition }x$$
leaves $dx$ invariant;  this property is stronger than that of incompressibility.  $\square$

 \end{example}

\vskip 2mm

In what follows, we shall denote by  ${\mathcal S}$  the set of incompressible semimartingales, 
by ${\mathcal S}_\nu$ the set of incompressible $\nu$-Brownian semimartingales and by ${\mathcal D}_\nu$ the set of incompressible $\nu$-Brownian
diffusions. Clearly we have

$${\mathcal D}_\nu \subset {\mathcal S}_\nu  \subset {\mathcal S}.$$

\begin{prop}\label{prop2.1} Let $g\in {\mathcal S}_\nu$, then for any $f\in C^2(M)$,
\begin{equation}\label{A2}
\E_{\mathbb P^g}(\langle\nabla f(g_t), u_t\rangle)=0.
\end{equation}
\end{prop}

\vskip 2mm
{\bf Proof.}  Taking the expectation with respect to ${\mathbb P^g}$ in \eqref{eq2.3}, we have
$$\E_{\mathbb P^g}(f(g_t))-\E_{\mathbb P^g}(f(g_0))
=\nu\int_0^t \E_{\mathbb P^g}(\Delta f(g_s))\, ds+\int_0^t\E_{\mathbb P^g}(\langle\nabla f(g_s), u_s\rangle)\,ds.$$
It follows that
$$\dis \nu \int_0^t \int_M \Delta f(x)\, dx\ ds+ \int_0^t\E_{\mathbb P^g}(\langle\nabla f(g_s), u_s\rangle)\,ds=0.$$
Since $\int_M \Delta f(x)\, dx=0$, we get the result. $\square$

\vskip 2mm
\begin{prop}\label{prop2.2} Let $g_t$ be a semimartingale on $M$ satisfying 
$$\dis dg_t(x)=\sum_{i=1}^m A_i(g_t(x))\circ dw_t^i + u_t(w,x)\, dx,$$
where $A_1, \cdots, A_m$ are $C^2$ divergence free vector fields on $M$ and $u_t(w,x)\in T_{g_t(x)}M$ is adapted such that $\int_M\E_x(\int_0^T|u_t(w,x)|^2dt)\,dx<+\infty$;
if $g$ is incompressible, then for any $f\in C^2(M)$
\begin{equation*}
\E_{\mathbb P^g}(\langle\nabla f(g_t), u_t\rangle)=0.
\end{equation*}

\end{prop}

\vskip 2mm
{\bf Proof.} Let $f\in C^2(M)$; then by It\^o formula \eqref{Ito}, 
\begin{equation*}
f(g_t)=f(g_0)+M^f_t + {1\over 2}\sum_{i=1}^m\int_0^t \Bigl( \langle \nabla_{A_i}(\nabla f), A_i\rangle 
+ \langle \nabla f, \nabla_{A_i}A_i\rangle \Bigr)ds+\int_0^t \langle\nabla f(g_s), u_s\rangle\, ds,
\end{equation*}

where $M_t^f$ is the martingale part. Note that $\dis \langle \nabla_{A_i}(\nabla f), A_i\rangle 
+ \langle \nabla f, \nabla_{A_i}A_i\rangle={\mathcal L}_{A_i}{\mathcal L}_{A_i}f$ where ${\mathcal L}_A$ denotes the Lie derivative with respect to $A$
; then taking the expectation under $\E_{\mathbb P}$, we get
$$\dis {1\over 2} \sum_{i=1}^m \Bigl(\int_M {\mathcal L}_{A_i}{\mathcal L}_{A_i}f\, dx \Bigr)+ \E_{\mathbb P^g}(\langle\nabla f(g_t), u_t\rangle)=0.$$
Since for each $i$, $\int_M {\mathcal L}_{A_i}{\mathcal L}_{A_i}f\, dx=0$, the result follows. $\square$

\vskip 2mm
In general  it is not clear whether the incompressibility condition implies the relation \eqref{A2}. However, the following is true:

\begin{prop}\label{prop3} Let $A_1, \cdots, A_m$ be $C^{2+\alpha} $ vector fields on $M$ and $A_0$ be a $C^{1+\alpha}$ vector field with some $\alpha>0$; consider
\begin{equation}\label{A3}
d\xi_t(x)=\sum_{i=1}^m A_i(\xi_t(x))\circ dw_t^i + A_0(\xi_t(x))\, dt,\quad \xi_0=x.
\end{equation}
Then for almost all $w$, the map $x\ra \xi_t(x)$ preserves the measure $dx$ if and only if $\div(A_i)=0$ for $i=0, 1, \cdots, m$.
\end{prop}

\vskip 2mm
{\bf Proof.}  We give a sketch of proof (see \cite{FangLuoTh} for more discussions). By \cite{Kunita}, 
; $x \ra \xi_t(x)$ is a diffeomorphism of $M$ and 
the push forward measure $(\xi_t^{-1})_\#(dx)$ of $dx$ by the inverse map of $\xi_t$ admits the density $K_t$ which is given by (see \cite{Kunita2}):

\begin{equation}\label{A4}
K_t(x)=\exp\Bigl( -\sum_{i=1}^m \int_0^t \div(A_i)(\xi_s(x))\circ dw_t^i -\int_0^t \div(A_0)(\xi_s(x))\, ds\Bigr).
\end{equation}

If $\div(A_i)=0$ for $i=0, 1, \cdots, m$, it is clear that $K_t=1$ and $x\ra \xi_t(x)$ preserves $dx$. Conversely, $K_t(x)=1$ for any $x\in M$ and $t\geq 0$  implies that,
$$\dis \sum_{i=1}^m \int_0^t \div(A_i)(\xi_s(x))\circ dw_t^i +\int_0^t \div(A_0)(\xi_s(x))\, ds=0;$$
or in It\^o form:
$$\dis \sum_{i=1}^m \int_0^t \div(A_i)(\xi_s(x)) dw_t^i +\int_0^t \Bigl[ {1\over 2}\sum_{i=1^m} {\mathcal L}_{A_i}\div(A_i)+ \div(A_0)\Bigr](\xi_s(x))\, ds=0.$$

The first term of above equality is of finite quadratic variation, while the second one is of finite variation; so that for each $i=1, \cdots, m$,
$\dis  \div(A_i)(\xi_s(x))=0$ and also 
$$\dis \Bigl[{1\over 2}\sum_{i=1^m} {\mathcal L}_{A_i}\div(A_i)+ \div(A_0)\Bigr](\xi_s(x))=0.$$

 It follows that, almost everywhere,
$$\dis \div(A_i)(\xi_s(x))=0 \quad\hbox{for}\quad i=0, 1, \cdots, m;$$
so that $\div(A_i)=0$ for $i=0, 1, \cdots, m$. $\square$

\vskip 2mm

According to \cite{Cruzeiro}, as well as \cite{AC2, Iyer, Eyink}, we introduce the following action functional  on semimartingales.

\begin{definition} Let

\begin{equation}\label{eq2.7}
S(g)=\frac{1}{2} \mathbb E_{\mathbb P^{g}}\left( \int_0^T |D_tg|^2 dt \right).
\end{equation}

\end{definition}

We say that $g$ has  {\it finite energy}
if $S(g)<\infty$. $\square$

\vskip 2mm

In what follows, we shall denote more precisely $D_tg(x)$ for $D_tg$ under the law $\mathbb P_x^g$. Then the action defined in \eqref{eq2.7} can be rewritten in the following form:

\begin{equation}\label{eq2.8}
S(g)= \frac{1}{2}\int_M \mathbb E_{\mathbb P^g_x}\left( \int_0^T |D_t g(x)|^2 dt \right) dx.
\end{equation}

We first recall briefly  known results about the calculus of stochastic  variation (see \cite{Cruzeiro, AC2, CCR}).
 Let $u_t(x)$ be a smooth vector field on a compact manifold (or on $\R^d$) which, for every $t$, is of divergence zero.
Consider an incompressible diffusion $g_t (x)$ with covariance $a$ such that $a (x,x)=2\mu \mathfrak{g}^{-1} (x)$ where $\mathfrak{g}$ is the metric tensor and time-dependent drift $u(t,\cdot )$. It defines a flow of diffeomorphisms preserving the volume measure.
We have: $D_tg(x)=u_t(g_t(x))$ and

$$\dis S(g)={1\over 2 }\int_{\T^d}\E_{\mathbb P_x^g}\Bigl(\int_0^T |u_t(g_t(x))|^2\,dt\Bigr)\,dx.$$

There are two manners to perform the perturbation. 
\vskip 2mm

{\it First perturbation of identity:}
\vskip 2mm

Let $w$ be a smooth  divergence free vector field  and $\alpha\in C^1(]0,T[)$. Consider, for for $\ee>0$, the ODE, 
\begin{equation}\label{eq2.10}
{d\Phi_t^\ee(x)\over dt}=\ee\,\alpha'(t)\, w(\Phi_t^\ee(x)), \Phi_0(x)=x.
\end{equation}

For each $t>0$, $\Phi_t^\ee$ is a perturbation of the identity map $id$. By It\^o's formula, for each fixed $\ee>0$, 
$t\ra \Phi_t^\ee(g_t(x))$ is a semimartingale starting from $x$.  Note that $g$ and $\Phi^\ee(g)$ are defined on the same probability space.
It was proved in \cite{Cruzeiro, AC2} that $u$ is a weak solution to Navier-Stokes equation if and only if $g$ is a
critical point of $S$. More precisely, $\dis {d\over d\ee}S(\Phi^\ee(g))_{|_{\ee=0}}=0$ if and only if
\begin{equation}\label{eq2.11}
\int_{\T^d}\!\!\int_0^T \langle u_t,  \alpha'(t) w+ \alpha(t)\, \nabla w\cdot u_t- \nu\, \alpha(t) \square w\rangle\, dtdx=0.
\end{equation}

\vskip 2mm

{\it Second perturbation of identity:}
\vskip 2mm
Note that in \cite{Iyer}, the perturbation of the identity was defined  in a different way. For each fixed $t>0$, the author of  \cite{Iyer} considered the ODE
\begin{equation}\label{perturbation}
{d\Psi^t_s\over ds}=\alpha(t) w(\Psi_s^t),\quad \Psi_0^t(x)=x.
\end{equation}
Set $\Psi(g)_t^\ee(x)=\Psi_\ee^t(g_t(x))$.
Then $\dis {d\over d\ee}S(\Psi(g)^\ee)_{|_{\ee=0}}=0$ if and only if the equation \eqref{eq2.11} holds.

Now we deal with the general case of compact Riemannian manifolds.

\begin{definition} Let  $M$ be a compact Riemannian manifold, $g$ a semimartingale on $M$ of finite energy. Define the probability measure $\mu$ on $[0,T]\times TM$ by
\begin{equation}\label{eq2.12}
\int_{[0,T]\times TM} f(t,x,v)\,\mu(dt,dx,dv)\\
={1\over T}\E_{\mathbb P^g}\Bigl[\int_0^T f\bigl(t, g(t), D_tg\bigr)dt\Bigr]
\end{equation}
where  $f: [0,T]\times TM\ra \R$ is any continuous function. $\square$
\end{definition}

We have the following result,

\begin{thm}\label{th1}
 Suppose that $ g\in {\mathcal S}_\nu$. Then  $g$ is a critical point of $ S$ with variations defined in \eqref{perturbation} if and ony if 
 $\mu$ is a solution to the Navier-Stokes equation in the sense of DiPerna-Majda,
 that is,
 
 \begin{equation}\label{eq2.13}
\int_0^T \int_{TM}\bigl[\alpha^{\prime} (t)\, v\cdot w  +\alpha (t)\, v\cdot\nabla_v w -\nu \alpha (t)\, 
v\cdot \hat\square w\bigr]\, d\mu (t,x,v)=0
\end{equation}
for all $\alpha \in C^1_c (]0,T[)$ and all smooth vector fields $w$ such that $\div(w)=0$. $\square$
\end{thm}

\vskip 2mm
{\bf Proof.} Let $\dis \Psi_\ee^t$ be the perturbation of identity defined in \eqref{perturbation}. Set $\eta_t^\ee=\Psi_\ee^t(g_t(x))$. Then $\{\eta_t^\ee,\ t\geq 0\}$ is 
a semimartingale on $M$. We denote by $(\xi^0(s), H_1(s), \ldots, H_m(s))$ the local characteristics of $g_t(x)$. 
By It\^o's formula (see \cite{Bismut}, p. 408), the drift term in local characteristics of $\eta_t^\ee$ is given by

\begin{equation}\label{eq2.14}
D_t  \Psi_\ee^t(g_t(x)) = {\partial \over \partial t}\Psi_\ee^t(g_t(x))+ d\Psi_\ee^t(g_t(x))\cdot \xi_t^0+{1\over 2}\sum_{i=1}^m \nabla(d\Psi_\ee^t)(g_t(x))(H_i(t),H_i(t)), 
\end{equation}
where $d\Psi_\ee^t(g_t(x))$ denotes the differential of $\Psi_\ee^t$ at $g_t(x)$. Let $\varphi(\ee,t)=D_t \Psi_\ee^t(g_t(x))\in T_{\eta_t^\ee}M $; then

\begin{equation*}
S(\Psi_\ee(g))={1\over 2} \E_{\mathbb P^g}\Bigl(\int_0^T \Bigl|\varphi(\ee,t)\Bigr|^2\, dt\Bigr).
\end{equation*}

We have: $\dis \varphi(0,t)=D_tg(x)$. Let 
$$\dis \varphi_1(\ee,t)= {\partial \over \partial t}\Psi_\ee^t(g_t(x)),$$
$$\dis \varphi_2(\ee,t)=d\Psi_\ee^t(g_t(x))\cdot \xi_t^0,$$
$$\dis \varphi_3(\ee,t)={1\over 2}\sum_{i=1}^m \nabla(d\Psi_\ee^t)(g_t(x))(H_i(t),H_i(t)).$$

Since the torsion is free, we have 
$$\dis {D\over d\ee}\varphi_1(\ee,t)_{|_{\ee=0}}={D\over d\ee}{d \over dt}\Psi_\ee^t(g_t(x))_{|_{\ee=0}}
={D\over dt}{d \over d\ee}_{|_{\ee=0}}\Psi_\ee^t(g_t(x))=\alpha'(t)w(x).$$

In order to compute the derivative of $\varphi_2$, consider a smooth curve $\beta(s)\in M$ such that $\beta(0)=g_t(x), \beta'(0)=D_tg(x)$. Then
$$\dis d\Psi_\ee^t(g_t(x))\cdot \xi_t^0={d\over ds}_{|_{s=0}}\Psi_\ee^t(\beta(s)).$$
Therefore

\begin{equation*}
\begin{split}
{D\over d\ee}_{|_{\ee=0}}\varphi_2(\ee,t)&={D\over ds}_{|_{s=0}}{d\over d\ee}_{|_{\ee=0}}\Psi_\ee^t(\beta(s))
={D\over ds}_{|_{s=0}}\Bigl[ \alpha(t)w(\beta(s))\Bigr]\\
&=\alpha(t)\,(\nabla w)(g_t(x))\cdot D_tg(x).
\end{split}
\end{equation*}

For computing $\varphi_3$, we shall use another description given in \cite{Bismut} (p. 405). For the moment, consider a $C^2$ map $f: M\ra M$.
 Let $x\in M$ and two tangent vectors $u, v\in T_xM$ be given. Let $x(t)\in M$ be a smooth curve such that $x(0)=x, x'(0)=u$, and $Y_t\in T_{x_t}M$ such that $Y_0=v$.
 Define $Q(f)(x): T_xM\times T_xM\ra T_{f(x)}M$ by
 
 \begin{equation}\label{eq2.16}
 Q(f)(x)(u,v)={d\over dt}_{|_{t=0}}\Bigl[ //_t^{-1} (df(x_t)\cdot Y_t)\Bigr]- df(x)\cdot \nabla_uv.
 \end{equation}
Then $\varphi_3$ can be expressed by

$$\dis \varphi_3(\ee,t)={1\over 2}\sum_{i=1}^m Q(\Psi_\ee^t(g_t(x)))(H_i(t), H_i(t)).$$

Let $\beta(s)\in M$ be a smooth curve such that $\beta(0)=g_t(x)$ and $\beta'(0)=H_i(t)$ and  $\{Y_s; s\geq 0\}$ be a family of tangent vectors along $\{\beta(s); s\geq 0\}$
such that $Y_0=H_i(t)$. Set

$$\dis \gamma(\ee, s)=\Psi_\ee^t(\beta(s))\quad\hbox{and} \quad X(\ee,s)=d\Psi_\ee^t(\beta(s))\cdot Y_s.$$

If $R$ denotes  be the curvature tensor on $M$, the following commutation relation holds,

\begin{equation*}
{D\over d\ee}{D\over ds}X(\ee,s)={D\over ds}{D\over d\ee}X(\ee,s)+ R\bigl({\partial \gamma\over \partial\ee}, {\partial \gamma\over \partial s}\bigr)X(\ee,s).
\end{equation*}

We have $\dis X(0,0)=H_i(t), {\partial \gamma\over \partial\ee}(0,0)=\alpha(t) w(x), {\partial \gamma\over \partial s}(0,0)=H_i(t)$;  therefore
$$\dis \Bigl[R\bigl({\partial \gamma\over \partial\ee}, {\partial \gamma\over \partial s}\bigr)X(\ee,s)\Bigr]_{|_{\ee=0, s=0}}=\alpha(t)\, R(w(g_t(x)), H_i(t))H_i(t).$$

Now let $c(\tau)\in M$ be a smooth curve such that $c(0)=\beta(s), c'(0)=Y_s$. We have
\begin{equation*}
\begin{split}
{D\over d\ee}_{|_{\ee=0}}X(\ee,s)&=\Bigl[{D\over d\tau}{d\over d\ee}\Psi_\ee^t(c(\tau))\Bigr](0,0)\\
&=\alpha(t) {D\over d\tau}_{|_{\tau=0}}w(c(\tau))=\alpha(t)\, (\nabla_{Y_s}w)(\beta(s)),
\end{split}
\end{equation*}

and
\begin{equation*}
{D\over ds}_{|_{s=0}}(\nabla_{Y_s}w)(\beta(s))=\langle \nabla_{H_i(t)}\nabla w, H_i(t)\rangle + \langle \nabla w, \nabla_{H_i(t)}H_i(t)\rangle.
\end{equation*}
Note that

$$\dis {D\over d\ee}_{|_{\ee=0}} d\Psi_\ee^t(g_t(x))\cdot \nabla_{H_i(t)}H_i(t)=\alpha(t)\,  \langle \nabla w, \nabla_{H_i(t)}H_i(t)\rangle.$$
 
Using \eqref{eq2.16}, we finally get
\begin{equation*}
 {D\over d\ee}_{|_{\ee=0}} \varphi_3(\ee,t)={1\over 2}\alpha(t)\, \sum_{i=1}^m \Bigl[ \langle\nabla_{H_i(t)}\nabla w, H_i(t)\rangle + R(w, H_i(t))H_i(t)\Bigr].
 \end{equation*}

When $g_t$ is a $\nu$-Brownian semimartingale, the right hand side of above equality is equal to
$$\dis \nu \alpha(t)\, (\Delta w + \Ric\, w )(g_t(x)),$$
which, due to \eqref{eq1.4-4}, is equal to 
$$\dis \nu\alpha(t)\, (-\hat\square w)(g_t(x)).$$

In conclusion,  $\dis {d\over d\ee}S(\Psi_\ee(g))_{|_{\ee=0}}=0$ yields

\begin{equation}\label{eq2.16bis}
\mathbb E_{\mathbb P^g} \int_0^T \Bigl[\alpha^{\prime} (t)\, w(g_t )\cdot D_t g +\alpha (t)\, 
(\nabla_{D_t g}w)(g_t )\cdot D_t g -\nu \alpha (t)\, \hat\square w (g_t )\cdot D_t g\Bigr]\,dt =0.
\end{equation}

According to \eqref{eq2.12}, the above equation  is nothing but \eqref{eq2.13}. $\square$

\vskip 2mm

As a consequence of this result, we obtain

\begin{thm}\label{th1bis} Let $(u_t)_{t\in [0,T]}$ be a family of divergence free  vector fields on $M$, which belong to the  Sobolev space ${\mathbb D}_1^2$ 
and are such that
\begin{equation}\label{eq2.16ter}
\int_M\!\!\int_0^T \Bigl( |u_t(x)|^2+ |\nabla u_t(x)|^2\Bigr)\, dtdx<+\infty;
\end{equation}
then equations \eqref{diffusion1}, \eqref{diffusion2} define an incompressible $\nu$-Brownian diffusion $\xi$ on $M$, which is a critical point of the action functional $S$
if and only if $u_t$ solves weakly the Navier-Stokes equation, that is,
\begin{equation}\label{Feb}
\int_{M}\!\!\int_0^T \langle u_t,  \alpha'(t) w+ \alpha(t)\, \nabla w\cdot u_t- \nu\, \alpha(t) \hat\square w\rangle\, dtdx=0
\end{equation}
\end{thm}
for all $\alpha \in C^1_c (]0,T[)$ and all smooth vector fields $w$ such that $\hbox{div}(w)=0$. $\square$

\vskip 2mm
{\bf Proof.} First we notice that in Proposition 4.3 in \cite{FangLuo}, the condition $q>2$  insures the tightness of a family
of probability measures; this condition can be relaxed to $q=2$
using Meyer-Zheng tightness results (see the proof of Theorem \ref{th2} below). Therefore by Theorem 6.4 in \cite{FangLuo}, 
 equations \eqref{diffusion1} and \eqref{diffusion2} define a diffusion process $\xi$, which is, a fortiori, in ${\mathcal S}_\nu$.
Therefore by the  above computations (see \eqref{eq2.16bis}), $\xi$ is a critical point to $S$ if and only if 

\begin{equation*}
\mathbb E_{\mathbb P^g} \int_0^T \Bigl[\alpha^{\prime} (t)\, w(\xi_t )\cdot u_t(\xi_t)  +\alpha (t)\, 
(\nabla_{u_t(\xi_t)}w)(\xi_t )\cdot u_t(\xi_t) -\nu \alpha (t)\, \hat\square w (\xi_t)\cdot u_t(\xi_t)\Bigr]\,dt =0,
\end{equation*}
which yields the result. \ $\square$

\begin{remark} It has been proved in \cite{Taylor} (see Theorem 4.6, p. 498) that for any $u_0\in L^2(M, dx)$, there exists 
$\{u_t, t\in [0,T]\}$ solution to \eqref{Feb}, satisfying Condition \eqref{eq2.16ter}. Therefore equations \eqref{diffusion1}, \eqref{diffusion2} 
define an incompressible $\nu$-Brownian diffusion $\xi$ on $M$, which is a critical point of the action functional $S$.
\end{remark}
\vskip 2mm

Note that in Theorem 3.2 of \cite{AC2}, a variational principle was established by using the first type of perturbations  of identity,  defined by \eqref{eq2.10}; 
on the other hand the manifold $M$ was supposed there to be a symmetric space in order to insure the existence of semimartingales with the desired properties.
A variational principe on a quite general Lie groups framework was derived in \cite{ACC} (c.f. also \cite{CCR}).
\vskip 2mm

In \cite{Brenier}, generalized flows  with prescribed initial and final configuration were introduced. It is quite difficult to construct incompressible semimartingales 
with    given prescriptions.
In order to emphasize the contrast with the situation in \cite{Brenier}, let's see the example of a Brownian bridge $g_t^{x,y}$ on $\R$ over $[0,1]$. 
 It is known that for $t<1$, $g_t^{x,y}$ solves the following SDE

\begin{equation}\label{B1}
dg_t^{x,y}=dw_t - {g_t^{x,y}-y\over 1-t}\, dt,\quad g_0^{x,y}=x.
\end{equation}

Then $g_t^{x,y}\ra y$ as $t\ra 1$ and we have
\begin{equation}\label{B2}
\E\Bigl(\int_0^1 |D_t g^{x,y}|^2\, dt\Bigr) =+\infty.
\end{equation}

Let $\eta$ be a probability measure on $M\times M$ having $dx$ as two marginals; we shall say that the incompressible semimartingale $\{g_t\}$ has $\eta$ as final configuration if
\begin{equation}\label{A5}
\E_{{\mathbb P}^g}(f(g_0, g_T))= \int_{M\times M} f(x,y)\, d\eta(x,y),\quad f\in C(M\times M).
\end{equation}

This means that the joint law of $(g_0, g_T)$ is $\eta$. If $g_t$  is as in Example \ref{ex2}, then
$$\dis \E_{{\mathbb P}^g}(f(g_0, g_T))=\int_{M\times M} f(x,y)p_T(x,y)\, dxdy,$$
where $p_t(x,y)$ is the heat kernel associated to $(g_t)$. Conversely if $(\rho_t(x,y))$ is solution to the following Fokker-Planck equation
\begin{equation*}
{d\over dt}\rho_t(x,y)=\nu\, \Delta_x \rho_t(x,y) + \langle u_t(x), \nabla_x\rho_t(x,y)\rangle,
\end{equation*}
with $\dis \lim_{t\ra 0}\rho_t=\delta_x$, for some $u\in L^2([0,T], {\mathbb D}_1^2(M))$ with $\div (u_t)=0$, we can construct an incompressible $\nu$-Brownian semimartingale
which has $\rho_T(x,y)dxdy$ as final configuration.

\vskip 2mm

In any case, we have the following result:

\begin{thm}\label{th2} Let $\eta$ be a probability measure as above.
If there exists an incompressible $\nu$-Brownian semimartingale $g$  on $M$ of finite energy $S(g)$ such that $\eta$ is its final configuration, then there exists one that minimizes 
the energy among all incompressible $\nu$-Brownian semimartingales having $\eta$ as final configuration.
\end{thm}

\vskip 2mm
{\bf Proof.}  Let $J: M\ra \R^N$ be an isometric embedding; then $dJ(x):\ T_xM\ra \R^N$ is such that for each $x\in M$ and $v\in T_xM$,
$\dis |dJ(x)\cdot v|_{\R^N} = |v|_{T_xM}$. Denote by $(dJ(x))^*: \R^N\ra T_xM$ the adjoint operator of $dJ(x)$, that is,
$$\dis \langle (dJ(x))^* a, v\rangle_{T_xM}=\langle dJ(x)v, a\rangle_{\R^N},\quad a\in\R^N, v\in T_xM.$$
Let $\{\ee_1, \ldots, \ee_N\}$ be an orthonormal basis of $\R^N$ and set
$$\dis A_i(x)=(dJ(x))^*\ee_i,\quad i=1, \ldots, N.$$

Then it is well-known that the vector fields $\{A_1, \ldots, A_N\}$ enjoy the following properties:

\vskip 2mm
\quad (i)  For any $v\in T_xM$, $\dis |v|_{T_xM}^2=\sum_{i=1}^N \langle A_i(x), v\rangle_{T_xM}^2$.

\quad (ii) $\dis \sum_{i=1}^N \nabla_{A_i}A_i=0$.

\vskip 2mm
Combining (i) and (ii) gives that $\dis \Delta_M f=\sum_{i=1}^N \L_{A_i}^2f$ for any $f\in C^2(M)$.
On the other hand, let $J(x)=(J_1(x), \ldots, J_N(x))$; then
$$\dis \langle dJ(x)v, \ee_i\rangle=dJ_i(x)\cdot v=\langle \nabla J_i(x), v\rangle_{T_xM},\quad\hbox{for any } v\in T_xM.$$

It follows that
\begin{equation}\label{S1}
A_i=\nabla J_i,\quad i=1, \cdots, N.
\end{equation}

Let $f\in C^2(M)$; then there exists $\bar f\in C^2(\R^N)$ such that $f(x)=\bar f(J(x))$. We have
\begin{equation}\label{S1.5}
\begin{split}
\L_{A_i}f&=\sum_{j=1}^N {\partial \bar f\over \partial x_j}(J(x))\,\langle \nabla J_j(x), A_i(x)\rangle\\
&=\sum_{j=1}^N {\partial \bar f\over \partial x_j}(J(x))\,\langle  A_j(x), A_i(x)\rangle.
\end{split}
\end{equation}
Therefore

\begin{equation*}
\begin{split}
\Delta_Mf=&\sum_{i=1}^N\sum_{j,k=1}^N {\partial^2\bar{f}\over \partial x_j\partial x_k}(J(x))\, \langle A_j, A_i\rangle\langle A_k, A_i\rangle\\
&+\sum_{i=1}^N\sum_{j=1}^N {\partial\bar{f}\over \partial x_j}(J(x))\L_{A_i}\langle A_j, A_i\rangle.
\end{split}
\end{equation*}

Notice that 
$$\dis \sum_{i=1}^N \L_{A_i}\langle A_j, A_i\rangle=\div(A_j)=\Delta_M J_j,$$
and according to property (i),
$$\dis \sum_{i=1}^N  \langle A_j, A_i\rangle\langle A_k, A_i\rangle=\langle A_j, A_k\rangle.$$
Finally the Laplacian $\Delta_M$ on $M$ can be expressed by

\begin{equation}\label{S2}
\Delta_M f=\sum_{j,k=1}^N {\partial^2\bar{f}\over \partial x_j\partial x_k}(J(x))\, \langle A_j, A_k\rangle 
+ \sum_{j=1}^N {\partial\bar{f}\over \partial x_j}(J(x))\, \Delta_M J_j.
\end{equation}

Having these preparations, we prove now the existence of a $g\in {\mathcal S}_\nu$ such that the minimum of action functinal $S$ 
is attained at $g$ in the class of those in ${\mathcal S}_\nu$ having $\eta$ as final configuration.
 Let 
$$\dis K=\inf_{g\in S_\nu} S(g).$$

There is a minimizing sequence $g^n\in {\mathcal S}_\nu$, that is,  $\dis \lim_{n\ra +\infty} S(g^n)=K$. Consider the canonical decomposition:

\begin{equation*}
J(g_t^n)= J(g_0^n) + M_t^n+ \int_0^t b^ n(s)\, ds.
\end{equation*}
Let $M_t^n=(M_t^{n,1}, \cdots, M_t^{n,N})$; then 

\begin{equation}\label{S3}
\langle M_t^{n,i}, M_t^{n,j}\rangle = 2\nu\ \int_0^t  \langle \nabla J_i, \nabla J_j\rangle(g_s^n)\,ds.
\end{equation}

By It\^o formula, we have
\begin{equation}\label{eq2.17}
b^n(t)= dJ(g_t^n) \cdot D_tg^n + \nu\, \Delta J(g_t^n).
\end{equation}

It follows that
\begin{equation*}
\E\Big(\int_0^T |b^n(t)|^2\, dt\Bigr)\leq 2 S(g^n) + 2T\nu\, ||\Delta J||_\infty.
\end{equation*}
Therefore $\dis \int_0^T |b^n(t)|^2\, dt$ is bounded in $L^2$. We can use Theorem 3 in \cite{Zheng} to conclude that the joint law $\hat P_n$ of 
$$\dis (J(g_\cdot^n), M_\cdot^n, B_\cdot^n, U_\cdot^n)$$
in $C([0,T], \R^N)\times C([0,T], \R^N)\times C([0,T], \R^N)\times C([0,T], \R^{N\times N})$ is a tight family, 
where 
$$\dis B_t^n=\int_0^t b^n(s)\,ds,\quad U_t^n=(\langle M_t^{n,i}, M_t^{n,j}\rangle )_{1\leq i,j\leq N}.$$
Let $\hat P$ be a limit point; up to a subsequence, we suppose that $\hat P_n$ converges weakly to $\hat P$. Again by Theorem 3 in \cite{Zheng}, under $\hat P$, the 
coordinate process
$$\dis (X_t, M_t, B_t, U_t)$$
has the following properties:

\vskip 2mm
\quad (i) $M_0=B_0=0, U_0=0$,

\quad (ii) $(M_t)$ is a local martingale such that $\dis U_t=(\langle M_t^{i}, M_t^{j}\rangle )_{1\leq i,j\leq N}$ and

\quad (iii)   $\dis B_t=\int_0^t b(s)\,ds$ with $\dis \int_0^T |b(s)|^2\, ds<+\infty$ almost surely.
\vskip 2mm

 Since $J(M)$ is closed in $\R^N$, we see that $X_t\in J(M)$. 
Let
$$\dis X_t=J(g_t).$$
For any $f\in C^2(M)$, by $\eqref{S2}$, we see that $f(g_t)$ is a real valued semimartingale. In other words, 
 $\{g_t;\ t\geq 0\}$ is a semimartingale on $M$.  Let $f\in C(M)$, the map $f\circ J^{-1}: J(M)\ra \R$ can be extended as a bounded continuous 
function on $\R^N$; therefore letting $n\ra \infty$, we get

\begin{equation*}
\int_M f(x)\, dx= \E(f(g^n(t)))=\E(f\circ J^{-1}(J(g^n_t)))\ra \E(f\circ J^{-1}(X_t))=\E(f(g_t)).
\end{equation*}

 In the same way, for $f\in C(M\times M)$,  we have
$$\dis \int_{M\times M} f(x,y)\,d\eta(x,y)=\E(f(g^n(0), g^n(T)))=\E(f(J^{-1}J(g^n(0)), J^{-1}J(g^n(T))))$$
which goes to, as $n\ra +\infty$,
$$\dis  \E(f(g(0), g(T))).$$
So $g$ is incompressible and has $\eta$ as final configuration. 

Besides, by \eqref{S3}, we have

 \begin{equation}\label{S3.5}
(\langle M_t^{i}, M_t^{j}\rangle )_{1\leq i,j\leq N}=2\nu\ \int_0^t  \langle \nabla J_i, \nabla J_j\rangle(g_s)\,ds.
\end{equation}

Let $f\in C^2(M)$; denote by $M_t^f$ the martingale part of $f(g_t)$. Then by It\^o formula,
$$\dis dM_t^f=\sum_{j=1}^N {\partial\bar{f}\over \partial x_j}(X_t)\, dM_t^j.$$
Therefore for $f_1, f_2\in C^2(M)$, according to \eqref{S3.5}, we have

$$\dis \langle dM^{f_1}_t,dM^{f_2}_t\rangle=\sum_{j,k=1}^N {\partial\bar{f_1}\over \partial x_j}(X_t){\partial\bar{f_2}\over \partial x_k}(X_t)
\, 2\nu \langle A_j, A_k\rangle_{g_t}\, dt.$$

On the other hand, using relation \eqref{S1.5} and property (i), we have

\begin{equation*}
\langle \nabla f_1, \nabla f_2\rangle=\sum_{\alpha=1}^N \L_{A_\alpha}f_1\, \L_{A_\alpha}f_2
=\sum_{j,k=1}^N  {\partial\bar{f_1}\over \partial x_j}{\partial\bar{f_2}\over \partial x_k}\langle A_j, A_k\rangle.
\end{equation*}
Combinant above two equalities, we finally get 

\begin{equation}\label{S3.6}
\langle dM^{f_1}_t,dM^{f_2}_t\rangle=2\nu\, \langle \nabla f_1, \nabla f_2\rangle_{g_t}\, dt.
\end{equation}

Since $X_t=J(g_t)$, we have 
$$\dis dB_t=dJ(g_t)\cdot D_tg\,dt + {1\over 2} \hbox{\rm Hess}J(g_t)\, dg_t\otimes dg_t.$$

Relation \eqref{S3.6} implies that ${1\over 2} \hbox{\rm Hess}J(g_t)\, dg_t\otimes dg_t=\nu \Delta_M J(g_t)\,dt$. Therefore we get
\begin{equation}\label{S4}
 B_t=\int_0^t dJ(g_s)\cdot D_sg\, ds + \nu \int_0^t \Delta_M J(g_s)\, ds.
\end{equation}
In conclusion $\{g_t;\ t\geq 0\}$ is a $\nu$-Brownian semimartingale on $M$ or
 $g\in {\mathcal S}_\nu$.

\vskip 2mm
We want to see that $K=S(g)$. Firstly
using the relation \eqref{eq2.17}, for any $t\in [0,T]$,
$$\dis \int_0^t dJ(g_s^n)\cdot D_sg^n\, ds = B_t^n-\nu \int_0^t \Delta J(g_s^n)\, ds.$$

Let $\phi: C([0,T],\R^N)\ra \R$ be a bounded continuous function, consider $\varphi: C([0,T],\R^N)\times C([0,T],\R^N)\ra \R$ defined by
$$\dis \varphi(B, g)=\phi\Bigl(B_\cdot-\nu\, \int_0^\cdot \Delta J(g_s)\,ds\Bigr).$$

Then $\varphi$ is a bounded continuous function on $C([0,T],\R^N)\times C([0,T],\R^N)$. It follows that $\dis \int_0^\cdot dJ(g_s^n)\cdot D_sg^n\, ds$ converges in law
to $\dis \int_0^\cdot dJ(g_s)\cdot D_sg\, ds$. 
 Let $\ee>0$;
 for $n$ big enough,

\begin{equation*}
\E\Big(\int_0^T |dJ(g_s^n)\cdot D_sg^n|^2\,ds\Bigr) \leq K+\ee.
\end{equation*}

Now by Theorem 10 in \cite{MeyerZheng}, 
\begin{equation*}
\E\Big(\int_0^T |dJ(g_s)\cdot D_sg|^2\,ds\Bigr) \leq K+\ee,
\end{equation*}
or $\dis \E\Big(\int_0^T | D_sg|^2\,ds\Bigr) \leq K+\ee$. Letting $\ee\ra 0$ gives $S(g)\leq K$. So $S(g)=K$.  $\square$



\vskip 2mm
\section{\bf   Classical solutions and generalized paths}\label{section3}

In this section, $M$ will be a torus: $M=\T^d$. Let $g\in {\mathcal D}_\nu$ be the solution of the following   SDE on $\T^d$
\begin{equation}\label{eq3.1}
dg_t=\sqrt{2\nu}\, dw_t - u(T-t, g_t)\,dt,\quad g_0\in \T^d
\end{equation}
where $g_0$ is a random variable having $dx$ as law, $w_t$ is the standard Brownian motion on $\R^d$, 
and $\{u(t,x); t\in [0,T]\}$ is a family of $C^2$ vector fields on $\T^d$, 
 identified to vector fields on $\R^d$ which are $2\pi$-periodic with respect to each space component. Suppose that $u$ is a 
strong solution to the Navier-Stokes equation
\begin{equation*}
{\partial\over \partial t}u(t,x) + \nabla u(t,x)\cdot u(t,x) -\nu\, \Delta u(t,x)=-\nabla p(T-t, x).
\end{equation*}

By It\^o 's formula, 
\begin{equation}\label{eq3.2}
\begin{split}
du(T-t, g_t)= &-({\partial u\over \partial t})(T-t, g_t)-\nabla u(T-t, g_t)\cdot u(T-t, g_t)\\&
 + \nu\, \Delta u(T-t,g_t)+\sqrt{2\nu}\, \nabla u(T-t, g_t)\cdot dw_t\\
&\hskip -3mm =\nabla p(t, g_t)\, dt + \sqrt{2\nu}\, \nabla u(T-t, g_t)\cdot dw_t
\end{split}
\end{equation}

According to definition \eqref{eq2.1}, $D_tg=-u(T-t,g_t)$ and 
\begin{equation}\label{eq3.3}
 D_tD_tg=-\nabla p(t, g_t).
\end{equation}

In what follows, we shall consider

\begin{equation}\label{eq3.4}
{\mathcal G}=\bigl\{g^*\in {\mathcal S}_\nu;\ dg_t^*=\sqrt{2\nu}\, dw_t+ D_tg^*\,dt,\ g^*(0)=g(0), g^*(T)= g(T)\bigr\}.
\end{equation}

Note that semimartingales in ${\mathcal G}$ are defined on a same probability space.

\begin{example}\label{ex3.1} Let $\alpha$ be a real continuous function on $\R^d$ and set
$$\dis \beta(w,t)=\sin({\pi t\over T})\int_0^t \alpha(w_s)\,ds,\quad c(w,t)={d\over dt}\beta(w,t).$$

Let $a\in \R^d$ be fixed. Consider $v(w,t)=c(w,t)a$; then $v$ is an adapted vector field on $\T^d$. Define
$$\dis g_t^*=g_t + \int_0^t v(w,s)\, ds.$$
Then $g^*\in {\mathcal G}$. $\square$
\end{example}

 We have the following result

\begin{thm}\label{th3.1}
Let  $g\in {\mathcal D}_\nu $ be  given in \eqref{eq3.1}.
Assume that the  process $g$ is associated with the Navier-Stokes equation in the sense that 
$$D_t D_t  g =-\nabla p(t, g_t )$$ 
a.s. for a {\it regular  pression} $p$ such that  $\nabla^2 p (t,x)\leq R\,\hbox{\rm Id}$, with $RT^2 \leq \pi^2$.
Then $g$ minimizes the energy $S$ in the class ${\mathcal G}$. $\square$
\end{thm}

\vskip 2mm

{\bf Proof. }
We define the following:

\begin{equation}\label{eq3.5}
B (g) =\frac{1}{2} \int_0^T |D_t g |^2 dt - \int_0^T p(t, g(t))dt
\end{equation}

Notice that the function $b(x,y)$ defined in \cite{Brenier} (p. 243) has no meaning in our setting (c.f.  \eqref{B1} and \eqref{B2}).
Let $g^*\in {\mathcal G}$; we shall prove that 

\begin{equation}\label{eq3.6}
\E(B(g))\leq \E(B(g^*)).
\end{equation}

Consider the function 
$$\dis \phi (t,x)=\frac{R}{2}|x|^2 -p(t,x).$$ 
For each $t\geq 0$, the function $x\ra \phi(t,x)$ is convex on $\R^d$ as $\nabla^2 p(t,x)\leq R\,\hbox{Id}$.
By It\^o formula

\begin{equation*}
d\bigl( D_tg\cdot g_t\bigr)=d(D_tg)\cdot g_t+\sqrt{2\nu}\, D_tg\cdot dw_t
+|D_tg|^2\, dt + d(D_tg)\cdot dg_t.
\end{equation*}

Analogously,
$$d(D_t g \cdot g_t^* )=d(D_t g)\cdot g_t^* +\sqrt{2\nu}\, D_t g \cdot dw_t+ D_tg\cdot D_tg^* +d(D_t g)\cdot dg_t^* .$$

Remarking that $\dis d(D_tg)\cdot dg_t= d(D_tg)\cdot dg_t^*$, and making the substraction of the above two equalities, we obtain

\begin{equation*}
d\bigl( D_tg\cdot (g_t^*-g_t)\bigr)= d(D_tg)\cdot (g_t^*-g_t)+ \bigl(D_tg\cdot D_tg^*-|D_tg|^2\bigr)\, dt.
\end{equation*}

It follows that
\begin{equation*}
\begin{split}
&D_Tg\cdot (g_T^*-g_T)- D_0g\cdot (g_0^*-g_0)\\
&\hskip -4mm=\int_0^T d(D_tg)\cdot (g_t^*-g_t) + \int_0^T \bigl(D_tg\cdot D_tg^*-|D_tg|^2\bigr)\, dt.
\end{split}
\end{equation*}

Notice that $g_0^*=g_0, g_T^*=g_T$, and using \eqref{eq3.1},  we have

\begin{equation}\label{eq3.7}
\begin{split}
&\int_0^T\bigl(-D_tg\cdot D_tg^*+|D_tg|^2\bigr)\, dt=\int_0^T d(D_tg)\cdot (g_t-g_t^*)\\
&=\int_0^T (g_t^*-g_t)\cdot \bigl( -\sqrt{2\nu}\, \nabla u(T-t, g_t)dw_t-\nabla p(t, g_t)dt\bigr).
\end{split}
\end{equation}

Using the  convexity, of $\phi$, we have

\begin{equation}\label{eq3.8}
\phi(t, g_t^*)-\phi(t, g_t)\geq \bigl( Rg_t-\nabla p(t,g_t)\bigr)\cdot (g_t^*-g_t).
\end{equation}

From \eqref{eq3.7} and \eqref{eq3.8}, we get
\begin{equation}\label{eq3.9}
\begin{split}
&\int_0^T\bigl(-D_tg\cdot D_tg^*+|D_tg|^2 + R g_t\cdot (g_t^*-g_t)\bigr)\, dt\\
&\hskip -5mm\leq -\sqrt{2\nu} \int_0^T (g_t^*-g_t)\cdot \nabla u(T-t, g_t)dw_t+\int_0^T \bigl(\phi(t,g_t^*)-\phi(t,g_t)\bigr)\, dt.
\end{split}
\end{equation}

We have
$\dis g_t^* -g_t=\int_0^t (D_s g^* -D_sg)\, ds$. Since $g^*_0-g_0=g_T^*-g_T=0$, by Poincar\'e 's inequaliy on the circle to get

$$\int_0^{T} |g_t^* -g_t|^2 dt \leq (\frac{T}{\pi})^2 \int_{0}^{T} |D_t g^* -D_t g|^2 dt.$$

Since $\dis ({T\over \pi})^2\leq {1\over R}$, we have
\begin{equation}\label{eq3.10}
{R\over 2}\int_0^{T} |g_t^* -g_t|^2 dt \leq {1\over 2}\int_{0}^{T} |D_t g^* -D_t g|^2 dt.
\end{equation}

Remark that the inequality, for $x,y, a, b\in\R$
$$\dis x^2-xy-Rb^2 + Rab\geq {1\over 2}x^2-{1\over 2}y^2 -{R\over 2}b^2 +{R\over 2} a^2$$
holds if and only if 
$$\dis {1\over 2} (x-y)^2 \geq {R\over 2} (b-a)^2.$$

Therefore by \eqref{eq3.10}, we have
\begin{equation}\label{eq3.11}
\begin{split}
&\int_0^T \bigl( |D_tg|^2 -D_tg\cdot D_tg^*-R|g_t|^2 + R g_t\cdot g_t^*\bigr)\,dt\\
&\geq \int_0^T \bigl( {1\over 2}|D_tg|^2 -{1\over 2}|D_tg^*|^2 -{R\over 2}|g_t|^2 +{R\over 2}|g_t^*|^2\bigr)\, dt.
\end{split}
\end{equation}

Combining \eqref{eq3.9} and \eqref{eq3.11}, we get
\begin{equation*}
\begin{split}
&\int_0^T \bigl( {1\over 2}|D_tg|^2 -{1\over 2}|D_tg^*|^2 -{R\over 2}|g_t|^2 +{R\over 2}|g_t^*|^2\bigr)\, dt\\
&\hskip -4mm  \leq -\sqrt{2\nu} \int_0^T (g_t^*-g_t)\cdot \nabla u(T-t, g_t)dw_t+\int_0^T \bigl(\phi(t,g_t^*)-\phi(t,g_t)\bigr)\, dt,
\end{split}
\end{equation*}
from which we deduce
\begin{equation*}
\begin{split}
&\int_0^T \bigl( {1\over 2}|D_tg|^2  -{R\over 2}|g_t|^2 +\phi(t,g_t)\bigr)\, dt\\
&\hskip -8mm  \leq -\sqrt{2\nu} \int_0^T (g_t^*-g_t)\cdot \nabla u(T-t, g_t)dw_t+\int_0^T \bigl({1\over 2}|D_tg^*|^2-{R\over 2}|g_t^*|^2+\phi(t,g_t^*)\bigr)\, dt,
\end{split}
\end{equation*}

or

\begin{equation*}
\begin{split}
&\int_0^T \bigl( {1\over 2}|D_tg|^2  -p(t,g_t)\bigr)\, dt\\
&\hskip -8mm  \leq -\sqrt{2\nu} \int_0^T (g_t^*-g_t)\cdot \nabla u(T-t, g_t)dw_t+\int_0^T \bigl({1\over 2}|D_tg^*|^2- p(t, g_t^*)\bigr)\, dt.
\end{split}
\end{equation*}

Using definition \eqref{eq3.5}, 
$$\dis B(g)\leq -\sqrt{2\nu} \int_0^T (g_t^*-g_t)\cdot \nabla u(T-t, g_t)dw_t+ B(g^*).$$
Taking the expectation of this inequality, we obtain \eqref{eq3.6}.  Notice that $\int_0^T \E(p(t,g_t))\, dt=\int_0^T \E(p(t, g_t^*))\,dt$; then \eqref{eq3.6} yields 
$\E(S(g))\leq \E(S(g^*))$. $\square$

\vskip 2mm
The following result provides a perturbation in a natural way and illustrates Theorem \ref{th3.1}.

\begin{prop} Let $v(w,t)$ be the vector field constructed in Example \ref{ex3.1}. Consider the following perturbation of $g_t$ given by \eqref{eq3.1}:
\begin{equation*}
dg_t^\ee = \sqrt{2\nu} dw_t -u(T-t, g_t)\, dt + \ee\, v(w,t)\, dt,\quad g_0^\ee=x.
\end{equation*}
Then we have
\begin{equation*}
{d\over d\ee} S(g^\ee)_{|_{\ee=0}}=0.
\end{equation*}
\end{prop}

\vskip 2mm
{\bf Proof.} We see that $\{g^\ee;\ \ee\geq 0\} \subset {\mathcal G}$. We have
\begin{equation*}
S(g^\ee)={1\over 2}\E\Bigl(\int_0^T |u(T-t,g_t) -\ee\, v(w,t)|^2\, dt \Bigr).
\end{equation*}

Therefore
\begin{equation*}
{d\over d\ee} S(g^\ee)_{|_{\ee=0}}=-\E\Bigl(\int_0^T\langle u(T-t, g_t), v(w,t)\rangle\, dt\Bigr).
\end{equation*}
Let $V_t=\int_0^t v_s\, ds$. By construction of $v$, $V_T=0$. Now by integration by parts,
\begin{equation*}
-\int_0^T \langle u(T-t, g_t), \dot V(w,t)\rangle\, dt = \int_0^T \langle d(u(T-t, g_t)), V(w,t)\rangle\, dt
\end{equation*}
which is equal to, using \eqref{eq3.2},
$ \int_0^T \langle \nabla p(t,g_t), V(w,t)\rangle\, dt$.
Therefore 
\begin{equation*}
\begin{split}
{d\over d\ee} S(g^\ee)_{|_{\ee=0}}=&\int_0^T \E\Bigl(\int_{\T^d} \langle \nabla p(t,g_t(x)), \beta(w,t) a\rangle\,dx\Bigr)\, dt\\
&=\int_0^T\E(\beta(w,t))\Bigl( \int_{\T^d} \langle \nabla p(t,x), a\rangle\, dx\Bigr)dt=0.
\end{split}
\end{equation*}
$\square$










\vskip 5mm

\begin{thebibliography}{W99}

\bibitem{AmbrosioF} L. Ambrosio, A. Figalli, Geodedics in the space of measure-preserving maps and plans, {\it Arch. Rational Mech. Anal.}, 194 (2009), 421--469.

\bibitem{AC1} A. Antoniouk, M. Arnaudon, A.B. Cruzeiro, Generalized stochastic flows and applications to incompressible viscous fluids, 
{\it Bull. Sci. Math.}, vol 138, issue 4 (2014), p. 565--584.

 \bibitem{ACC} M. Arnaudon, X. Chen and A.B. Cruzeiro, Stochastic Euler-Poincar\'e reduction,
{\it J. Math. Physics}, 55 (2014), 081507.

\bibitem{AC2} M. Arnaudon, A.B. Cruzeiro, Lagrangian Navier-Stokes diffusions on manifolds: variational principle and stability, 
{\it Bull. Sci. Math.}, 136 (8) (2012), p. 857--881.

\bibitem{AC3} M. Arnaudon, A.B. Cruzeiro, Stochastic Lagrangian flows on some compact manifolds, {\it Stochastics}, 84 (2-3) (2012), p. 367--381.

\bibitem{Arnold} V. I. Arnold, Sur la g\'eom\'etrie diff\'erentielle des groupes de Lie de dimension infinie et ses applications
 \`a l' hydrodynamique des fluides parfaits, {\it Ann. Inst. Fourier}, 16 (1966), 316--361.


\bibitem{Brenier} Y. Brenier, The least action principle and the related concept of generalized flows for incompressible perfect fluids.
{\it J. AMS}, 2 (1989), 225-255.

\bibitem{Brenier2} Y. Brenier, Minimal geodesics on groups of volume-preserving maps and generalized solutions of the Euler equations, 
{\it Comm. Pure and Appl. Maths}, 52 (1999), 411--452.

\bibitem{Bismut} J.M. Bismut, M\'ecanique al\'eatoire, {\it Lect. Notes in Maths}, 866, 1981

\bibitem{CCR} X. Chen, A.B. Cruzeiro, T. Ratiu, Constrained and stochastic variational principles for dissipative equations with advected quantities,
arXiv:1506.05024 

\bibitem{Constantin} P.Constantin, G. Iyer, A stochastic Lagrangian representation of the three-dimensional incompressible Navier-Stokes equations,
{\it Comm. Pure Appl. Math.}, 61 (2008), 330-345.


\bibitem{Cruzeiro} F. Cipriano and A.B. Cruzeiro, { Navier-Stokes equations and diffusions on the group of homeomorphisms
of the torus}, {\it Comm. Math. Phys.}  275 (2007), 255--269.

\bibitem{DiPerna} R. DiPerna, A. Majda, Oscillations and concentrations in weak solutions of incompressible fluid equations, 
{\it Comm. Math. Phys.}, 108 (1987), 667--689.

\bibitem{EM} D.G. Ebin, J.E. Marsden, Groups of diffeomorphisms and the motion of an incompressible fluid,
{\it Ann. of Math.} 92 (1970), 102-163.

\bibitem{Eyink} G. L. Eyink, Stochastic least-action principle for the incompressible Navier-Stokes equation, {\it Physica D}, 239 (2010), 1236--1240.

\bibitem{FangLuo}  S. Fang, H.Li, D. Luo, Heat semi-group and generalized flows on complete Riemannian manifolds. {\it Bull. Sci. Math.},  135 (2011), 565--600.

\bibitem{FangLuoTh} S. Fang, D. Luo, A. Thalmaier,  Stochastic differential equations with coefficients in Sobolev spaces. {\it J. Funct. Analysis}, 259 (2010), 1129--1168.

\bibitem{Fleming} W.H. Fleming, H.M. Soner, Controlled Markov processes and viscosity solutions, Stoch. Modelling and Applied Prob. 25,  \emph{Springer, Berlin}, 2006.
(2nd ed.)

\bibitem{Follmer} H. Follmer, Random fields and diffusion processes.
\emph{\'Ecole d'\'et\'e de Probabilit\'es de Saint-Flour, XV--XVII--1982--87}, 143--303, Lecture Notes in Math., 1362, \emph{Springer, Berlin}, 1988.

\bibitem{Iyer} G. Iyer, A variational principle for the Navier-Stokes equations, preprint, 2008.

\bibitem{Kunita} H. Kunita, Stochastic differential equations and stochastic flows of diffeomorphisms. 
\emph{\'Ecole d'\'et\'e de Probabilit\'es de Saint-Flour, XII--1982}, 143--303, Lecture Notes in Math., 1097, \emph{Springer, Berlin}, 1984.

\bibitem{Kunita2} H. Kunita, {\it Stochastic Flows and Stochastic
Differentail Equations}. Cambridge University Press, 1990.

\bibitem{Leonard} C. L\'eonard, A survey of the Schr\"odinger problem and some of its connections with optimal transport, \emph{Discrete Contin. Dyn. Syst.}, 34 (2014), 
1533--1574.

\bibitem{MeyerZheng}
P.A. Meyer, W.A. Zheng,  {Tightness criteria for laws of semimartingales}, {\it Annales de l'I.H.P.}, section~B, tome~20, n.~4 (1984), 353--372.

\bibitem{MT} M. Mitrea, M. Taylor, Navier-Stokes equations on Lipschitz domains in Riemannian manifolds, {\it Math. Ann.}, 321(2001), 955-987.

\bibitem{Pier} V. Pierfelice, The incompressible Navier-Stokes equations on non-compact manifolds, arXive: 1406.1644v1.

\bibitem{Taylor} M. Taylor,  \emph{Partial Differential Equations III: Nonlinear Equations, Nonlinear equations}, 
Vol. 117, Applied Mathematical Sciences, Springer New York second edition (2011).

\bibitem{Yasue}
K. Yasue, A variational principle for the Navier-Stokes equation, \textit{J. Funct. Anal.},  51 (2), (1983) 133--141.

\bibitem{Zambrini} J.C. Zambrini, Variational processes and stochastic versions of Mechanics. {\it J. Math. Phys} 27 (9), (1986), 2307

\bibitem{Zheng}  W.~A.~Zheng,
Tightness results for laws of diffusion processes. {\it Ann. de l' I.H.P.}, section B {21} (2) (1985), 103--124.

\end{thebibliography}
\end{document}